\documentclass{amsart}

\usepackage{amsmath,amsthm,amsfonts,amscd,amssymb}
\usepackage[pdfborder={0 0 0}]{hyperref}
\usepackage{sagetex}

\newtheorem{thm}{Theorem}[section]
\newtheorem{cor}[thm]{Corollary}
\newtheorem{lem}[thm]{Lemma}

\newtheorem{quest}{Question}

\theoremstyle{definition}

\theoremstyle{remark}
\newtheorem{rem}{Remark}[section]

\begin{document}

\title{Bounds for solid angles of lattices of rank three}
\author{Lenny Fukshansky \and Sinai Robins}

\address{Department of Mathematics, Claremont McKenna College, 850 Columbia Avenue, Claremont, CA 91711-6420}
\email{lenny@cmc.edu}
\address{Division of Mathematical Sciences, Nanyang Technological University, SPMS-MAS-03-01, 21 Nanyang Link, Singapore 637371 }
\email{rsinai@ntu.edu.sg}
\subjclass[2000]{Primary: 52C07, 05B30; Secondary: 11H06, 51M25}
\keywords{lattices, solid angles, spherical configurations}

\begin{abstract} We find sharp absolute constants $C_1$ and $C_2$ with the following property: every well-rounded lattice of rank 3 in a Euclidean space has a minimal basis so that the solid angle spanned by these basis vectors lies in the interval $[C_1,C_2]$. In fact, we show that these absolute bounds hold for a larger class of lattices than just well-rounded, and the upper bound holds for all. We state a technical condition on the lattice that may prevent it from satisfying the absolute lower bound on the solid angle, in which case we derive a lower bound in terms of the ratios of successive minima of the lattice. We use this result to show that among all spherical triangles on the unit sphere in $\mathbb R^N$ with  vertices on the minimal vectors of a lattice, the smallest possible area is achieved by a configuration of minimal vectors of the (normalized) face centered cubic lattice in $\mathbb R^3$. Such spherical configurations come up in connection with the kissing number problem.
\end{abstract} 

\maketitle

\def\A{{\mathcal A}}
\def\B{{\mathcal B}}
\def\C{{\mathcal C}}
\def\D{{\mathcal D}}
\def\E{{\mathcal E}}
\def\F{{\mathcal F}}
\def\x{{\mathcal H}}
\def\I{{\mathcal I}}
\def\J{{\mathcal J}}
\def\K{{\mathcal K}}
\def\L{{\mathcal L}}
\def\Ll{{\mathfrak L}}
\def\M{{\mathcal M}}
\def\Mm{{\mathfrak M}}
\def\Pp{{\mathfrak P}}
\def\Aa{{\mathfrak A}}
\def\Ss{{\mathfrak S}}
\def\N{{\mathcal N}}
\def\PP{{\mathcal P}}
\def\R{{\mathcal R}}
\def\s{{\mathcal S}}
\def\V{{\mathcal V}}
\def\W{{\mathcal W}}
\def\X{{\mathcal X}}
\def\Y{{\mathcal Y}}
\def\H{{\mathcal H}}
\def\OO{{\mathcal O}}
\def\aaa{{\mathbb A}}
\def\cee{{\mathbb C}}
\def\Nn{{\mathbb N}}
\def\pee{{\mathbb P}}
\def\que{{\mathbb Q}}
\def\real{{\mathbb R}}
\def\zed{{\mathbb Z}}
\def\gmn{{\mathbb G_m^N}}
\def\qbar{{\overline{\mathbb Q}}}
\def\DL{{\underline{\Delta}}}
\def\DU{{\overline{\Delta}}}
\def\eps{{\varepsilon}}
\def\vek{{\varepsilon_k}}
\def\ahat{{\hat \alpha}}
\def\bhat{{\hat \beta}}
\def\gt{{\tilde \gamma}}
\def\h{{\tfrac12}}
\def\ba{{\boldsymbol a}}
\def\be{{\boldsymbol e}}
\def\bei{{\boldsymbol e_i}}
\def\bc{{\boldsymbol c}}
\def\bm{{\boldsymbol m}}
\def\bk{{\boldsymbol k}}
\def\bi{{\boldsymbol i}}
\def\bl{{\boldsymbol l}}
\def\bq{{\boldsymbol q}}
\def\bu{{\boldsymbol u}}
\def\bt{{\boldsymbol t}}
\def\bs{{\boldsymbol s}}
\def\bv{{\boldsymbol v}}
\def\bw{{\boldsymbol w}}
\def\bx{{\boldsymbol x}}
\def\bX{{\boldsymbol X}}
\def\bz{{\boldsymbol z}}
\def\bwy{{\boldsymbol y}}
\def\bg{{\boldsymbol g}}
\def\bY{{\boldsymbol Y}}
\def\bL{{\boldsymbol L}}
\def\baa{{\boldsymbol\alpha}}
\def\bb{{\boldsymbol\beta}}
\def\bet{{\boldsymbol\eta}}
\def\bxi{{\boldsymbol\xi}}
\def\bo{{\boldsymbol 0}}
\def\bol{{\boldsymbol 1}_L}
\def\ep{\varepsilon}
\def\p{\boldsymbol\varphi}
\def\q{\boldsymbol\psi}
\def\WR{\operatorname{WR}}
\def\rank{\operatorname{rank}}
\def\aut{\operatorname{Aut}}
\def\lcm{\operatorname{lcm}}
\def\sgn{\operatorname{sgn}}
\def\spn{\operatorname{span}}
\def\md{\operatorname{mod}}
\def\Norm{\operatorname{Norm}}
\def\dim{\operatorname{dim}}
\def\det{\operatorname{det}}
\def\Vol{\operatorname{Vol}}
\def\rk{\operatorname{rk}}
\def\md{\operatorname{mod}}
\def\sqp{\operatorname{sqp}}
\def\Aut{\operatorname{Aut}}
\def\GL{\operatorname{GL}}
\def\Sim{\operatorname{Sim}}

\section{Introduction}
\label{intro}

Given an integer $N \geq 2$, let $\real^N$ be the Euclidean space with the usual norm $\|\ \|$ on it. The kissing number problem in dimension $N$ asks for the maximal number of non-overlapping unit balls in $\real^N$ that touch another unit ball. The answer to this problem is currently only known in dimensions $N=2,3,4,8,24$ (see \cite{bachoc}, \cite{conway}, \cite{musin}, and \cite{odlyzko}). In fact, in dimension 3 this was the subject of a famous argument between Isaac Newton and David Gregory, where Newton claimed that the kissing number is 12 while Gregory believed it was 13; two different proofs that this number is 12 finally appeared in the 1950-s, by Sch\"{u}tte and van der Waerden \cite{waerden} and by Leech \cite{leech} (see also \cite{conway}, \cite{hsiang} for details, including an extensive bibliography), although there were previous unsuccessful and incomplete attempts. The kissing number problem can be reformulated as follows: find the maximal configuration of points on the unit sphere in $\real^N$ such that the {\it angular separation} between any pair of these points is at least $\pi/3$; by angular separation between two points we mean the smallest angle between the vectors connecting these two points to the origin in the plane spanned by these vectors. Such configurations are often expected to come from sets of {\it minimal vectors} (i.e. vectors of smallest nonzero norm) of lattices, at least this is the case in all dimensions where the kissing number is known. For example, in dimension three an optimal configuration of 12 points is given by the set of minimal vectors of the (normalized) face centered cubic lattice $A_3$. Define a {\it spherical lattice-minimal triangle} to be a non-degenerate spherical triangle with vertices at the endpoints of minimal vectors of a lattice; as we will discuss below, the minimality condition forces the angular separation between every pair of these points to be at least $\pi/3$. The above consideration raises the following naturally related question.

\begin{quest} \label{quest1} Given a spherical lattice-minimal triangle on a unit sphere in $\real^N$, what is the minimal possible two-dimensional spherical area it can have?
\end{quest}

\noindent
In other words, although minimal vectors of the normalized lattice $A_3$ form the largest (with respect to cardinality) configuration of points on the unit sphere with angular separation at least $\pi/3$, and hence produce spherical triangles of small area (at least on the average), could it be that minimal vectors of some other lattice form some spherical triangle of even smaller area? Our first result answers this question.

\begin{thm} \label{sp_area} Any spherical lattice-minimal triangle on a unit sphere in $\real^N$ has area at least 
\begin{equation}
\label{sp1}
0.551285598\hdots = 4\pi \times 0.043869914\hdots,
\end{equation}
which is precisely the area of the spherical triangle formed by the vectors
\begin{equation}
\label{sp_sharp}
\left( \begin{matrix} 1/\sqrt{2} \\ 1/\sqrt{2} \\ 0 \end{matrix} \right),\ \left( \begin{matrix} 0 \\ 1/\sqrt{2} \\ 1/\sqrt{2} \end{matrix} \right),\ \left( \begin{matrix} 1/\sqrt{2} \\ 0 \\ 1/\sqrt{2} \end{matrix} \right)
\end{equation}
on the unit sphere in $\real^3$. These are precisely minimal basis vectors of the face centered cubic lattice $A_3$, normalized to lie on the unit sphere.
\end{thm}

To prove Theorem \ref{sp_area}, we use a somewhat more general and technical result, which we present next. We start with some basic notation. Let $B_N$ be a unit ball in $\real^N$. Given a lattice $\Lambda \subset \real^N$ of rank~$r$, we define its successive minima
$$0 < \lambda_1 \leq \dots \leq \lambda_r$$
by
$$\lambda_i = \inf \{ \lambda \in \real_{>0} : \Lambda \cap \lambda B_N \text{ contains } i \text{ linearly independent vectors} \}.$$
There exists a collection of linearly independent vectors $\bx_1,\dots,\bx_r$ in $\Lambda$ such that $\|\bx_i\| = \lambda_i$ for each $1 \leq i \leq r$; we will refer to them as vectors {\it corresponding} to successive minima. When $r \leq 4$, these vectors form a basis for $\Lambda$, which is precisely a minimal basis; this is not necessarily true for $r \geq 5$ (see for instance \cite{pohst}). Let us also write
$$S(\Lambda) = \left\{ \bx \in \Lambda : \|\bx\| = \lambda_1 \right\}$$
for the set of all minimal vectors of $\Lambda$. In the special case when $S(\Lambda)$ contains $r$ linearly independent vectors, i.e. when
$$\lambda_1 = \dots = \lambda_r,$$
$\Lambda$ is called a {\it well-rounded} lattice, abbreviated WR. WR lattices and configurations of their minimal vectors play an important role in discrete optimization problems (see \cite{martinet}, \cite{conway}, \cite{hsiang}). In particular, spherical configurations which give good kissing numbers always come from WR lattices.

We will say that a basis $\bx_1, \dots, \bx_r$ for a lattice $\Lambda \subset \real^N$ of rank $r$ is {\it minimal} (often referred to as {\it Minkowski reduced} in the literature - see for instance \cite{lek}) if $\bx_1 \in S(\Lambda)$ and for each $2 \leq i \leq r$, $\bx_i \in \Lambda$ is a shortest vector such that the collection of vectors $\bx_1,\dots,\bx_i$ is extendable to a basis in~$\Lambda$. Minimal bases for lattices are extensively studied for their importance in number theory and discrete geometry, as well as in computer science and engineering applications (see \cite{pohst}, \cite{near:ort}, \cite{esm}). In particular, if a lattice happens to have an orthogonal basis, it is easy to prove that this basis must be minimal. In fact, it turns out that even suitably defined ``near-orthogonality" of the basis vectors is sufficient for a basis to be minimal (or ``almost minimal" - here we mean reduced vs Minkowski reduced). This idea has been successfully exploited by the famous LLL algorithm \cite{LLL}. It should be mentioned that in addition to the notion of Minkowski reduction that we are using here other basis reduction procedures exist, most notably the Hermite-Korkine-Zolotareff reduction (see for instance \cite{martinet} for details); the common general principle behind various reduction procedures is the minimization of the orthogonality defect.

In the reverse direction, one can ask how ``close" to orthogonal does a minimal basis have to be? When $r=2$, the answer is classical: if $\bx_1,\bx_2$ is a minimal basis for a lattice $\Lambda \subset \real^N$ of rank two, then the (normalized) angle between these vectors has to lie in the interval~$[1/6, 1/3]$; moreover, for every such lattice there exists a minimal basis with the (normalized) angle between vectors lying in the interval~$[1/6, 1/4]$ (since we will use this result, we include its proof in section~\ref{2D} to make our presentation self-contained). It is natural to ask for analogues of this statement for lattices of higher rank, starting with~$r=3$ as follows.

\begin{quest} \label{quest2} Do there exist absolute constants $C_1$ and $C_2$ such that every lattice of rank three has a minimal basis that spans a solid angle in the interval $[C_1,C_2]$?
\end{quest}

Notice that when $\Lambda$ is a WR lattice of rank three in $\real^N$, then $S(\Lambda)$ contains a minimal basis. In fact, by scaling if necessary, we can assume that $\lambda_1 = \lambda_2 = \lambda_3 = 1$, and then any three linearly independent vectors in $S(\Lambda)$ form a minimal basis and span a spherical lattice-minimal triangle on the unit sphere centered at the origin in $\real^N$. On the other hand, vertices of a spherical lattice-minimal triangle on the unit sphere centered at the origin in $\real^N$ form a minimal basis for a WR lattice of rank three with $\lambda_1 = \lambda_2 = \lambda_3 = 1$. Since the solid angle in question is equal to the normalized spherical area of the corresponding spherical triangle, there is a direct connection between Question \ref{quest2} restricted to just WR lattices and Question \ref{quest1}. In this note, we answer Question \ref{quest2} completely for a wide class of lattices including all WR lattices and provide a partial answer (with the constant $C_1$ depending on the lattice) for all the remaining lattices. We then use this result, along with our method, to answer Question \ref{quest1}.

Let $\Lambda$ be a lattice of rank 3 in $\real^N$ and let $\bx_1,\bx_2,\bx_3$ be a fixed minimal basis for $\Lambda$, corresponding to successive minima $\lambda_1,\lambda_2,\lambda_3$, respectively. We will write $\theta_{ij}$ for the angle between the vectors $\bx_i$ and $\bx_j$, $1 \leq i < j \leq 3$, and will refer to $\theta_{12},\theta_{13},\theta_{23}$ as the {\it vertex angles} corresponding to this minimal basis. We also define the ratios of successive minima
\begin{equation}
\label{KKK}
K_{12} = \frac{\lambda_1}{\lambda_2},\ K_{23} =  \frac{\lambda_2}{\lambda_3},\ K_{13} = \frac{\lambda_1}{\lambda_3} = K_{12} K_{23}.
\end{equation}
Then clearly $\max \{ K_{12}, K_{23} \} \leq 1$, and $\min \{ K_{12}, K_{23} \} < 1$ if and only if $\Lambda$ is not WR. Let us also define
\begin{equation}
\label{nu}
\nu := \frac{1}{4} \cos^{-1} \left( \frac{1}{2} \left( 2K_{12} K_{23} \cos \theta_{12} + 2K_{12} \cos \theta_{13} - K_{12}^2 K_{23} - K_{23}  \right) \right) < \frac{\pi}{6}.
\end{equation}

\begin{rem} \label{solid_measure} Our measure of three-dimensional solid angles is always normalized to be between 0 and 1, i.e. we divide by $4\pi$, the maximal possible solid angle which corresponds to the full sphere in $\real^3$.
\end{rem}

\noindent
We can now state our next result.

\begin{thm} \label{main_3} Let $\Lambda \subset \real^N$ be a lattice of rank 3. Then there exists a minimal basis $\bx_1,\bx_2,\bx_3$ for $\Lambda$ such that
\begin{equation}
\label{mn3.1}
\pi/3 \leq \theta_{12}, \theta_{13} \leq \pi/2,\ \pi/3 \leq \theta_{23} \leq 2\pi/3.
\end{equation}
Fix any minimal basis $\bx_1,\bx_2,\bx_3$ for $\Lambda$ satisfying (\ref{mn3.1}), and let $\Omega$ be the solid angle  formed by these vectors. Then $\Omega \leq 0.125$. If in addition
\begin{equation}
\label{mn3.2}
\theta_{23} \leq  \cos^{-1} \left( \cos \theta_{12} + \cos \theta_{13} - 1 \right),
\end{equation}
then $\Omega \geq 0.043869914\hdots$ The above condition is satisfied in particular by every WR lattice, and the bounds on $\Omega$ are sharp as demonstrated in Remark \ref{sharp} below. On the other hand, if
$$\cos^{-1} \left( \cos \theta_{12} + \cos \theta_{13} - 1 \right) < \theta_{23} \leq 2\pi/3,$$
then in fact
$$\theta_{23} \leq \cos^{-1} \left( \frac{1}{2} \left( 2K_{12} K_{23} \cos \theta_{12} + 2K_{12} \cos \theta_{13} - K_{12}^2 K_{23} - K_{23}  \right) \right) < \frac{2\pi}{3},$$
and
$$\Omega \geq \frac{1}{\pi} \tan^{-1} \left( \tan \nu\ \sqrt{\frac{1 - 3 \tan^2 \nu}{3 - \tan^2 \nu}} \right) > 0,$$
where $\nu$ is as in \eqref{nu}.
\end{thm}

In section~\ref{2D} we discuss some basic two-dimensional lemmas, which will be later used in our main argument. Theorem \ref{main_3} is proved in section~\ref{3D}: it follows immediately by combining Lemmas \ref{min_bas}, \ref{solid_3}, \ref{solid_other}, Corollary \ref{solid_wr}, and Remark \ref{sharp}. In  section~\ref{sphere}, we use the techniques developed in section~\ref{3D} to prove Theorem~\ref{sp_area}. We are now ready to proceed.
\bigskip

\section{Preliminary two-dimensional lemmas}
\label{2D}

We start with some basic lemmas about the angles between minimal basis vectors for lattices of rank~2, which we will use in our argument for lattices of rank~3.

\begin{lem} \label{vec_angles} Let $\bx_1$ and $\bx_2$ be nonzero vectors in $\real^N$ so that the angle $\theta$ between them satisfies either $0 < \theta < \frac{\pi}{3}$ or $\frac{2\pi}{3} < \theta < \pi$. Then
$$\min \left\{ \| \bx_1 - \bx_2 \|,  \| \bx_1 + \bx_2 \| \right\} < \max \{ \| \bx_1 \|, \| \bx_2 \| \}.$$
\end{lem}

\proof
Notice that either
$$\frac{1}{2} < \cos \theta = \frac{\bx_1^t \bx_2}{\|\bx_1\| \|\bx_2\|} < 1,\text{ or } -1 < \cos \theta = \frac{\bx_1^t \bx_2}{\|\bx_1\| \|\bx_2\|} < -\frac{1}{2}.$$
In the first case:
\begin{eqnarray*}
\|\bx_1 - \bx_2\|^2 & = & (\bx_1-\bx_2)^t (\bx_1-\bx_2) = \|\bx_1\|^2 + \|\bx_2\|^2 - 2 \bx_1^t \bx_2 \\
& < & \|\bx_1\|^2 + \|\bx_2\|^2 - \|\bx_1\| \|\bx_2\| < \max \{ \| \bx_1 \|, \| \bx_2 \| \}^2.
\end{eqnarray*}
In the second case:
\begin{eqnarray*}
\|\bx_1 + \bx_2\|^2 & = & (\bx_1+\bx_2)^t (\bx_1+\bx_2) = \|\bx_1\|^2 + \|\bx_2\|^2 + 2 \bx_1^t \bx_2 \\
& < & \|\bx_1\|^2 + \|\bx_2\|^2 - \|\bx_1\| \|\bx_2\| < \max \{ \| \bx_1 \|, \| \bx_2 \| \}^2.
\end{eqnarray*}
\endproof

\begin{lem} \label{angle} Let $\Lambda \subset \real^N$ be a lattice of rank 2 with successive minima $\lambda_1 \leq \lambda_2$, and let $\bx_1,\bx_2$ be vectors in $\Lambda$ corresponding to $\lambda_1,\lambda_2$, respectively. Let $\theta$ be the angle between $\bx_1$ and $\bx_2$. Then
$$\pi/3 \leq \theta \leq 2\pi/3.$$
\end{lem}

\proof
Clearly $\theta \in (0, \pi)$. Assume that either $0 < \theta < \pi/3$ or $\frac{2\pi}{3} < \theta < \pi$, then Lemma \ref{vec_angles} implies that
$$\min \left\{ \| \bx_1 - \bx_2 \|,  \| \bx_1 + \bx_2 \| \right\} < \| \bx_2 \| = \lambda_2,$$
which contradicts the definition of $\lambda_2$ since the vectors $\bx_1$ and $\bx_1\pm \bx_2$ are linearly independent.
\endproof

\begin{lem} \label{plane_angles} Let $\Lambda \subset \real^N$ be a lattice of full rank with successive minima 
$$0 < \lambda_1 \leq \dots \leq \lambda_N,$$
and vectors $\bx_1,\dots,\bx_N$ corresponding to these successive minima, respectively, chosen in such a way that all of them lie in the half-space $x_N \geq 0$. For every pair of indices $1 \leq i < j \leq N$, let $\theta_{ij}$ be the angle between the vectors $x_i$ and $x_j$. Then
$$\pi/3 \leq \theta_{ij} \leq 2 \pi/3.$$
\end{lem}

\proof
Let $\Lambda_{ij} = \spn_{\zed} \{ \bx_i, \bx_j \}$, then the successive minima $\mu_1, \mu_2$ of $\Lambda_{ij}$ are:
$$\mu_1 = \lambda_i \leq \mu_2 = \lambda_j.$$
Indeed, it is clear that $\mu_1 \leq \mu_2$ and $\mu_1 \leq \lambda_i$, $\mu_2 \leq \lambda_j$, so suppose for instance that $\mu_1 < \lambda_i$. Then there exists a vector $\bwy = a_i \bx_i + a_j \bx_j$, such that $0 \neq a_i, a_j \in \zed$ and
\begin{equation}
\label{a1}
\|\bwy\| = \mu_1 < \lambda_i.
\end{equation}
Then the collection of vectors
\begin{equation}
\label{a2}
\bx_1,\dots,\bx_{i-1},\bwy,\bx_{i+1},\dots,\bx_N
\end{equation}
must be linearly independent in $\Lambda$, and so by definition of successive minima, $\mu_1 \geq \lambda_i$, which is a contradiction. If, on the other hand, we assume that $\mu_2 < \lambda_j$, then we can apply the same argument as above replacing $\mu_1$ with $\mu_2$ in (\ref{a1}) and inserting $\bwy$ for $\bx_j$ instead of $\bx_i$ in (\ref{a2}) to reach the same contradiction. Therefore $\lambda_i,\lambda_j$ are the successive minima of $\Lambda_{ij}$, and hence $\bx_i,\bx_j$ are vectors corresponding to successive minima. Now the conclusion follows by Lemma \ref{angle} above.
\endproof

\bigskip

\section{Bounds for solid angles: proof of Theorem \ref{main_3}}
\label{3D}

In this section we prove a collection of lemmas, which together comprise the result of Theorem \ref{main_3}. First we need to fix our choice of a minimal basis for a lattice of rank~3, which is accomplished by the following lemma.

\begin{lem} \label{min_bas}  Let $\Lambda \subset \real^N$ be a lattice of rank 3. Then there exists a minimal basis $\bx_1,\bx_2,\bx_3$ for $\Lambda$ such that
\begin{equation}
\label{min_angle}
\pi/3 \leq \theta_{12}, \theta_{13} \leq \pi/2,\ \pi/3 \leq \theta_{23} \leq 2\pi/3,
\end{equation}
where $\theta_{ij}$ is the angle between vectors $\bx_i$ and $\bx_j$ for all $1 \leq i < j \leq 3$.
\end{lem}

\proof
Let 
$$0 < \lambda_1 \leq \lambda_2 \leq \lambda_3$$
be successive minima of $\Lambda$, and let $\bx_1,\bx_2,\bx_3$ be the vectors corresponding to $\lambda_1, \lambda_2, \lambda_3$ respectively. Then Lemma \ref{plane_angles} implies that $\pi/3 \leq \theta_{ij} \leq 2\pi/3$. In fact, if $\theta_{12} > \pi/2$, replace $\bx_2$ with $- \bx_2$, and if $\theta_{13} > \pi/2$, replace $\bx_3$ with $- \bx_3$. Then $\bx_1,\bx_2,\bx_3$ is a minimal basis satisfying (\ref{min_angle}), as required.
\endproof

In what follows, the measure of three-dimensional solid angles is normalized as specified in Remark~\ref{solid_measure}, however we assume it is converted back to steradians when we compute values of trigonometric functions of $\Omega$.

\begin{lem} \label{solid_3} Let $\Lambda \subset \real^N$ be a lattice of rank 3, and let $\bx_1,\bx_2,\bx_3$ be a minimal basis for $\Lambda$ guaranteed by Lemma \ref{min_bas}. Let $\Omega$ be the solid angle  formed by these vectors. Then $\Omega \leq 0.125$. If in addition
\begin{equation}
\label{cos-1}
\theta_{23} \leq  \cos^{-1} \left( \cos \theta_{12} + \cos \theta_{13} - 1 \right),
\end{equation}
then $\Omega \geq 0.043869914\hdots$
\end{lem}

\proof
With our choice of the basis $\bx_1,\bx_2,\bx_3$, we have
\begin{equation}
\label{tij.0}
\pi/3 \leq \theta_{12}, \theta_{13} \leq \pi/2,\ \pi/3 \leq \theta_{23} \leq 2\pi/3,
\end{equation}
Now make the choice of $\pm \bx_3$ that makes $\theta_{23}$ as small as possible for (\ref{tij.0}) to hold; notice in particular that if $\theta_{13} = \pi/2$, then $\pm \bx_3$ can be chosen to ensure that $\pi/3 \leq \theta_{23} \leq \pi/2$. From now on we will always assume this choice of $\bx_1, \bx_2, \bx_3$. Let
\begin{equation}
\label{ts1}
\alpha = \frac{\theta_{12} + \theta_{13}}{4},\ \beta = \frac{\theta_{12}-\theta_{13}}{4},\ c = \frac{\theta_{23}}{4},
\end{equation}
and so
\begin{equation}
\label{ts2.0}
\frac{\pi}{6} \leq \alpha \leq \frac{\pi}{4},\ -\frac{\pi}{24} \leq \beta \leq \frac{\pi}{24},\ \frac{\pi}{12} \leq c \leq \frac{\pi}{6}.
\end{equation}
Let $\Omega$ be the solid angle spanned by the vectors $\bx_1,\bx_2,\bx_3$, then L'Huilier's theorem (see, for instance, \cite{sphere_trig}) implies that
\begin{equation}
\label{tan1.1}
\tan \left( \frac{\Omega}{4} \right)^2 = \tan \left( \alpha + c \right) \tan \left( \alpha - c \right) \tan \left( c+ \beta \right) \tan \left( c - \beta \right),
\end{equation}
when $\Omega$ is measured in steradians. Notice that for each $c \in (-\pi/4,\pi/4)$,
\begin{equation}
\label{beta}
\tan \left( c + \beta \right) \tan \left( c - \beta \right) = \frac{\tan^2 c - \tan^2 \beta}{1 - \tan^2 c \tan^2 \beta}
\end{equation}
is a decreasing function of $|\beta|$, and so (\ref{beta}) is minimized when $\beta$ is as large as possible, and maximized when $\beta = 0$.
\smallskip

To produce an upper bound for $\Omega$, let us write
\begin{equation}
\label{xy}
x_1 = \pi/2 - \theta_{12},\ x_2 = \pi/2 - \theta_{13},\ y = \theta_{23} - \pi/2, \text{ and } x = x_1+x_2,
\end{equation}
so that 
\begin{equation}
\label{xy1}
\alpha = \frac{\pi}{4} - \frac{x}{4},\ c = \frac{\pi}{8} + \frac{y}{4}, \text{ where } 0 \leq x \leq \frac{\pi}{3},\ -\frac{\pi}{6} \leq y \leq \frac{\pi}{6},
\end{equation}
and $y$ is as small as possible for (\ref{xy1}) to hold. Then (\ref{tan1.1}), along with the fact that (\ref{beta}) is maximized when $\beta=0$, implies that
\begin{equation}
\label{tan1.3.2}
\tan \left( \frac{\Omega}{4} \right)^2 \leq \tan \left( \frac{3\pi}{8} - \frac{x}{4} + \frac{y}{4} \right) \tan \left( \frac{\pi}{8} - \frac{x}{4} - \frac{y}{4}  \right) \tan \left( \frac{\pi}{8} + \frac{y}{4} \right)^2
\end{equation}
It is not difficult to observe that with the constraints of (\ref{xy1}) satisfied, the right hand side of (\ref{tan1.3.2}) is a decreasing function of $x$, so to maximize it we can assume that $x=0$, meaning that $\theta_{12} = \theta_{13} = \pi/2$. In this case $\pi/3 \leq \theta_{23} \leq \pi/2$, meaning that $-\pi/6  \leq y \leq 0$, and
\begin{equation}
\label{tan1.3.3}
\tan \left( \frac{\Omega}{4} \right)^2 \leq \tan \left( \frac{3\pi}{8} + \frac{y}{4} \right) \tan \left( \frac{\pi}{8} - \frac{y}{4}  \right) \tan \left( \frac{\pi}{8} + \frac{y}{4} \right)^2,
\end{equation}
where the right hand side of (\ref{tan1.3.3}) is an increasing function of $y$. Therefore
$$\tan \left( \frac{\Omega}{4} \right)^2 \leq \tan \left( \frac{3\pi}{8} \right) \tan \left( \frac{\pi}{8} \right)^3,$$
and the upper bound for $\Omega$ follows.
\smallskip

Next we produce a lower bound for $\Omega$. Assume that in addition to (\ref{tij.0}), the inequality (\ref{cos-1}) is also satisfied. Then
\begin{equation}
\label{tij}
\pi/3 \leq \theta_{12}, \theta_{13} \leq \pi/2,\ \pi/3 \leq \theta_{23} \leq  \min \left\{ \cos^{-1} \left( \cos \theta_{12} + \cos \theta_{13} - 1 \right), 2\pi/3 \right\},
\end{equation}
and $\bx_1,\bx_2,\bx_3$ is chosen so that $\theta_{23}$ is as small as possible for (\ref{tij}) to hold. Then
\begin{equation}
\label{ts2}
\frac{\pi}{6} \leq \alpha \leq \frac{\pi}{4},\ -\frac{\pi}{24} \leq \beta \leq \frac{\pi}{24},\ \frac{\pi}{12} \leq c \leq \min \left\{ \frac{1}{4} \cos^{-1} \left( 2 \cos 2\alpha \cos 2 \beta - 1 \right), \frac{\pi}{6} \right\}.
\end{equation}
The right hand side of (\ref{tan1.1}) is easily checked to be an increasing function of $c$ when the constraints of (\ref{ts2}) are satisfied (we use Sage mathematical software package \cite{sage} - see illustration in section~\ref{sage_code}). Therefore for each fixed pair of values of $\alpha$ and $\beta$,
\begin{equation}
\label{tan1.2}
\tan \left( \alpha + \frac{\pi}{12} \right) \tan \left( \alpha - \frac{\pi}{12} \right) \tan \left( \frac{\pi}{12}+ \beta \right) \tan \left( \frac{\pi}{12} - \beta \right) \leq \tan \left( \frac{\Omega}{4} \right)^2,
\end{equation}
when (\ref{cos-1}) holds. Since $\theta_{12}, \theta_{13} \geq \pi/3$, it follows that
$$4 |\beta| = |\theta_{12}-\theta_{13}| \leq 4 \alpha - 2\pi/3,$$
and so (\ref{beta}) is minimized when $\beta = \min \{ \alpha - \pi/6, \pi/24\}$. We can also assume without loss of generality that $0 < \Omega < 1/2$ (i.e. the measure of $\Omega$ in steradians is between 0 and $2\pi$), and so $\tan \left( \frac{\Omega}{4} \right)^2$ is an increasing function of $\Omega$.

First assume that $\alpha - \pi/6 \leq \pi/24$, then $\alpha \leq 5\pi/24$ and the above observations combined with (\ref{tan1.2}) imply that
\begin{equation}
\label{tan1.3}
\tan \left( \frac{\Omega}{4} \right)^2 \geq \tan \left( \alpha + \frac{\pi}{12} \right) \tan \left( \alpha - \frac{\pi}{12} \right)^2 \tan \left( \frac{\pi}{4} - \alpha \right) \geq \tan \left( \frac{\pi}{12} \right)^3.
\end{equation}

Next assume that $\alpha - \pi/6 \geq \pi/24$, then $\alpha \geq 5\pi/24$ and the above observations combined with (\ref{tan1.2}) imply that
\begin{eqnarray}
\label{tan1.3.1}
\tan \left( \frac{\Omega}{4} \right)^2 & \geq & \tan \left( \alpha + \frac{\pi}{12} \right) \tan \left( \alpha - \frac{\pi}{12} \right) \tan \left( \frac{\pi}{8} \right) \tan \left( \frac{\pi}{24} \right) \nonumber \\
& \geq & \tan \left( \frac{7\pi}{24} \right) \tan \left( \frac{\pi}{8} \right)^2 \tan \left( \frac{\pi}{24} \right) > \tan \left( \frac{\pi}{12} \right)^3.
\end{eqnarray}
The lower bound for $\Omega$ now follows in the case when (\ref{cos-1}) is satisfied.
\endproof
\smallskip

\begin{cor} \label{solid_wr} Let $\Lambda \subset \real^N$ be a well-rounded lattice of rank 3, and let $\bx_1,\bx_2,\bx_3$ be a minimal basis for $\Lambda$ guaranteed by Lemma \ref{min_bas}. Then (\ref{cos-1}) is satisfied, and so
$$0.043869914\hdots \leq \Omega \leq 0.125,$$
by Lemma \ref{solid_3}.
\end{cor}

\proof
Let $\lambda$ be the common value of successive minima of $\Lambda$, and notice that
$$\| \bx_1-\bx_2-\bx_3 \|^2 = \lambda^2(3 + 2(-\cos \theta_{12} - \cos \theta_{13} + \cos \theta_{23} )) \geq \lambda^2,$$
which means that $-\cos \theta_{12} - \cos \theta_{13} + \cos \theta_{23}  \geq -1$.
\endproof

\begin{rem} \label{sharp} The bounds of Lemma \ref{solid_3} are sharp. The lower bound is achieved by the fcc (face centered cubic) lattice $A_3$ (same as $D_3$) in $\real^3$ with the choice of a minimal basis
$$\bx_1 = \left( \begin{matrix} 1 \\ 1 \\ 0 \end{matrix} \right),\ \bx_2 = \left( \begin{matrix} 0 \\ 1 \\ 1 \end{matrix} \right),\ \bx_3 = \left( \begin{matrix} 1 \\ 0 \\ 1 \end{matrix} \right),$$
so that the angles between these vectors are
$$\theta_{12} = \theta_{13} = \theta_{23} = \pi/3.$$
The upper bound is achieved by the integer lattice $\zed^3$ in $\real^3$ with the standard choice of a minimal basis
$$\bx_1 = \left( \begin{matrix} 1 \\ 0 \\ 0 \end{matrix} \right),\ \bx_2 = \left( \begin{matrix} 0 \\ 1 \\ 0 \end{matrix} \right),\ \bx_3 = \left( \begin{matrix} 0 \\ 0 \\ 1 \end{matrix} \right),$$
so that the angles between these vectors are
$$\theta_{12} = \theta_{13} = \theta_{23} = \pi/2.$$
Notice that both of these lattices are well-rounded, and so in fact the bounds of Corollary \ref{solid_wr} are sharp.
\end{rem}
\bigskip

Next we consider the situation when (\ref{cos-1}) is not satisfied, in which case we can derive a lower bound for the solid angle formed by the minimal basis vectors of a lattice depending on the ratios of successive minima.

\begin{lem} \label{solid_other} Let $\Lambda \subset \real^N$ be a lattice of rank 3, let $\bx_1,\bx_2,\bx_3$ be a minimal basis for $\Lambda$, corresponding to successive minima $\lambda_1, \lambda_2, \lambda_3$, as guaranteed by Lemma \ref{min_bas}, and let $\Omega$ be the solid angle  formed by these vectors. By (\ref{tij.0}), we can assume that
$$\pi/3 \leq \theta_{12}, \theta_{13} \leq \pi/2,\ \pi/3 \leq \theta_{23} \leq 2\pi/3,$$
and $\theta_{23}$ is as small as possible for these inequalities to hold. Suppose in addition that
\begin{equation}
\label{cos-1-2}
\cos^{-1} \left( \cos \theta_{12} + \cos \theta_{13} - 1 \right) < \theta_{23} \leq 2\pi/3.
\end{equation}
Let $K_{12}$, $K_{13}$, and $K_{23}$ be as in \eqref{KKK}, then $\min \{ K_{12}, K_{23} \} < 1$ and
\begin{equation}
\label{O-nu}
\Omega \geq \frac{1}{\pi} \tan^{-1} \left( \tan \nu\ \sqrt{\frac{1 - 3 \tan^2 \nu}{3 - \tan^2 \nu}} \right) > 0,
\end{equation}
where
\begin{equation}
\label{nu_def}
\nu = \frac{1}{4} \cos^{-1} \left( \frac{1}{2} \left( 2K_{12} K_{23} \cos \theta_{12} + 2K_{12} \cos \theta_{13} - K_{12}^2 K_{23} - K_{23}  \right) \right) < \frac{\pi}{6}.
\end{equation}
\end{lem}

\proof
We use the notation from the proof of Lemma \ref{solid_3}. The first implication of (\ref{cos-1-2}) is that $-\frac{1}{2} \leq \cos \theta_{23}  < -1 + \cos \theta_{12} + \cos \theta_{13}$, so in particular $\cos \theta_{12} + \cos \theta_{13} \geq 1/2$. Combining this observation with (\ref{tij.0}), we obtain
\begin{equation}
\label{tij_2}
\frac{\pi}{3} \leq \theta_{12} \leq \frac{\pi}{2},\ \frac{\pi}{3} \leq \theta_{13} \leq \cos^{-1} \left( \frac{1}{2} - \cos \theta_{12} \right),
\end{equation}
in addition to (\ref{cos-1-2}). 

Another implication of (\ref{cos-1-2}) is that $\Lambda$ cannot be well-rounded, by Corollary \ref{solid_wr}. Therefore $\min \{ K_{12}, K_{23} \} < 1$, and thus $K_{13} < 1$. Notice that the vectors $\bx_1,\bx_2,\bx_1-\bx_2-\bx_3$ are linearly independent, and hence
\begin{eqnarray*}
\lambda_3^2 & \leq & \|\bx_1-\bx_2-\bx_3\|^2 \\ 
& = & \lambda_3^2 \left( 1 + K_{13}^2 + K_{23}^2 - 2K_{13} K_{23} \cos \theta_{12} - 2K_{13} \cos \theta_{13} + 2K_{23} \cos \theta_{23} \right).
\end{eqnarray*}
Therefore, since $\cos \theta_{13} \geq 1/2 - \cos \theta_{12}$, $\cos \theta_{12} \leq 1/2$, and $K_{13} = K_{12} K_{23}$,
\begin{eqnarray*}
\cos \theta_{23} & \geq & -\frac{1}{2K_{23}} \left( K_{13}^2 + K_{23}^2 - 2K_{13} K_{23} \cos \theta_{12} - 2K_{13} \cos \theta_{13} \right) \\
& = & -\frac{1}{2} \left( K_{12}^2 K_{23} + K_{23} - 2K_{12} K_{23} \cos \theta_{12} - 2K_{12} \cos \theta_{13} \right) \\
& \geq & -\frac{1}{2} \left( K_{12}^2 K_{23} + K_{23} - K_{12} + 2K_{12} (1 - K_{23}) \cos \theta_{12}  \right) \\
& \geq & -\frac{K_{23}}{2} \left( K_{12}^2  -  K_{12} + 1  \right) > - \frac{1}{2},
\end{eqnarray*}
since $\min \{ K_{12}, K_{23} \} < 1$ and $\max \{ K_{12}, K_{23} \} \leq 1$. Therefore
\begin{equation}
\label{cos_K}
\theta_{23} \leq \cos^{-1} \left( \frac{1}{2} \left( 2K_{12} K_{23} \cos \theta_{12} + 2K_{12} \cos \theta_{13} - K_{12}^2 K_{23} - K_{23}  \right) \right) < \frac{2\pi}{3}.
\end{equation}
\smallskip

Let now $\nu$ be as in (\ref{nu_def}). It can now be easily checked that the right hand side of (\ref{tan1.1}) is a decreasing function of $c$ when the constraints of (\ref{cos-1-2}) and (\ref{tij_2}) are satisfied (we use Sage mathematical software package \cite{sage} - see illustration in section~\ref{sage_code}), and $c \leq \nu < \pi/6$ by (\ref{cos_K}). Therefore for each fixed pair of values of $\alpha$ and $\beta$,
\begin{equation}
\label{tan1.4}
\tan \left( \alpha + \nu \right) \tan \left( \alpha - \nu \right) \tan \left( \nu + \beta \right) \tan \left( \nu - \beta \right) \leq \tan \left( \frac{\Omega}{4} \right)^2,
\end{equation}
when (\ref{cos-1-2}) holds. The left hand side of (\ref{tan1.4}) is now minimized when $\alpha$ is as small as possible. Notice that
$$4\alpha \geq 2\pi/3 + 4|\beta|,$$
and hence
\begin{equation}
\label{tan1.5}
 \tan \left( \frac{\Omega}{4} \right)^2 \geq \tan \left( \frac{\pi}{6} + |\beta| + \nu \right) \tan \left( \frac{\pi}{6} + |\beta| - \nu \right) \tan \left( \nu + \beta \right) \tan \left( \nu - \beta \right),
\end{equation}
where $0 \leq |\beta| \leq \pi/24$. The right hand side of (\ref{tan1.5}) is easily checked to be an increasing function of $|\beta|$, so to minimize take $|\beta| = 0$, therefore
\begin{equation}
\label{tan1.6}
\tan \left( \frac{\Omega}{4} \right)^2 \geq \left( \frac{1 - 3 \tan^2 \nu}{3 - \tan^2 \nu} \right) \tan^2 \nu,
\end{equation}
and (\ref{O-nu}) now follows from (\ref{tan1.6}), since  the left hand side of (\ref{tan1.6}) is an increasing function of $\Omega$, as was indicated in the proof of Lemma \ref{solid_3}.
\endproof
\bigskip

\section{Area of spherical triangles: proof of Theorem~\ref{sp_area}}
\label{sphere}

In this section we prove Theorem~\ref{sp_area}. Consider a spherical lattice-minimal triangle $T$ on the unit sphere centered at the origin in $\real^N$, and let $\bx_1,\bx_2,\bx_3$ be minimal vectors of some lattice $L$ corresponding to the vertices of this triangle. Let $\Omega$ be the three-dimensional solid angle spanned by $\bx_1,\bx_2,\bx_3$, measured in steradians (i.e. not normalized as in Remark~\ref{solid_measure}), then the measure of $\Omega$ is precisely the spherical two-dimensional area of $T$. Hence we want to show that $\Omega$ is greater or equal than the number in (\ref{sp1}). 

Let $\Lambda = \spn_{\zed} \{ \bx_1,\bx_2,\bx_3 \}$. Then 
$\bx_1,\bx_2,\bx_3$ are minimal vectors in $\Lambda$, since they are minimal vectors in $L$ and $\Lambda \subseteq L$, and
$$\|\bx_1\| = \|\bx_2\| = \|\bx_3\| = 1.$$
Since these vectors form a spherical triangle, they must be linearly independent, and so they form a minimal basis for $\Lambda$ with solid angle $\Omega$, in particular $\Lambda$ is WR with all the successive minima equal to 1. If the vertex angles $\theta_{12},\theta_{13},\theta_{23}$ corresponding to these vectors satisfy~(\ref{mn3.1}), then Theorem~\ref{main_3} implies the result with equality precisely in the case described in (\ref{sp_sharp}), which is a minimal basis for the normalized fcc lattice (see Remark~\ref{sharp}).

Now assume that the vertex angles do not satisfy~(\ref{mn3.1}). This means that at least two of the angles $\theta_{12}, \theta_{13}, \theta_{23}$ must lie in the interval $(\pi/2,2\pi/3]$. Since $\Lambda$ is WR we can reindex the three vectors as necessary, and so we can assume without loss of generality that
$$\pi/2 < \theta_{12} \leq \theta_{23} \leq 2\pi/3,\ \pi/3 \leq \theta_{13} \leq \theta_{12}.$$
Let $\alpha$, $\beta$, and $c$ be as in (\ref{ts1}), then
\begin{equation}
\label{sp_bounds}
\frac{5\pi}{24} < \alpha \leq \frac{\pi}{3},\ 0 \leq \beta \leq \frac{\pi}{12},\ \frac{\pi}{8} \leq c \leq \frac{\pi}{6}.
\end{equation}
Now $\Omega$ is given by (\ref{tan1.1}), and since it is proportional to the spherical area of our triangle, we want to understand how small can it be. The right hand side of (\ref{tan1.1}) is easily checked to be an increasing function of $c$ when the constraints of (\ref{sp_bounds}) are satisfied (we use Sage mathematical software package \cite{sage} - see illustration in section~\ref{sage_code}); also, as discussed in the proof of Lemma~\ref{solid_3}, to minimize $\Omega$ we need to maximize $|\beta|$, so we can take $\beta=\pi/12$. Therefore for each fixed value of $\alpha$ and $\beta$,
\begin{eqnarray}
\label{sp_tan1}
\tan \left( \frac{\Omega}{4} \right)^2 & \geq & \tan \left( \alpha + \frac{\pi}{8} \right) \tan \left( \alpha - \frac{\pi}{8} \right) \tan \left( \frac{\pi}{8}+ \beta \right) \tan \left( \frac{\pi}{8} - \beta \right) \nonumber \\
& \geq &  \tan \left( \alpha + \frac{\pi}{8} \right) \tan \left( \alpha - \frac{\pi}{8} \right) \tan \left( \frac{5\pi}{24} \right) \tan \left( \frac{\pi}{24} \right) \nonumber \\
& > & \tan \left( \frac{\pi}{3} \right) \tan \left( \frac{5\pi}{24} \right) \tan \left( \frac{\pi}{12} \right) \tan \left( \frac{\pi}{24} \right) > \tan \left( \frac{\pi}{12} \right)^3,
\end{eqnarray}
which means that $\Omega$ is greater than the number in (\ref{sp1}), and hence finishes the proof of the theorem.
\bigskip

\section{Appendix: Sage code}
\label{sage_code}

Here we present the Sage code illustrating the increasing/decreasing behavior of functions used in the arguments above. Here is the code corresponding to the proof of Lemma~\ref{solid_3}:

\begin{sageblock}
#auto
@interact
def plotter(x = slider(pi/3,pi/2), y = slider(pi/3,pi/2)):
        z = var('z')
        a, b, c = (x+y)/4, (x-y)/4, z/4
        expression(z) = tan(a+c)*tan(a-c)*tan(c+b)*tan(c-b)
        zmin, zmax = pi/3, min(2*pi/3,arccos(cos(x)+cos(y)-1))
        color = (1,.25,0)
        plot(expression, z, zmin, zmax, rgbcolor=color)
            .show(xmin=pi/2,xmax=2*pi/3,ymin=0,ymax=0.3)
\end{sageblock}

\noindent
Here is the code corresponding to the proof of Lemma~\ref{solid_other}:

\begin{sageblock}
#auto
@interact
def plotter(x = slider(pi/3,pi/2), y = slider(pi/3,pi/2)):
    if arccos(-1+cos(x)+cos(y)) <= 2*pi/3:        
        z = var('z')
        a, b, c = (x+y)/4, (x-y)/4, z/4
        expression(z) = tan(a+c)*tan(a-c)*tan(c+b)*tan(c-b)
        zmin, zmax = arccos(-1+cos(x)+cos(y)), 2*pi/3
        color = (1,.25,0)
        plot(expression, z, zmin, zmax, rgbcolor=color)
            .show(xmin=pi/3,xmax=2*pi/3,ymin=0,ymax=.1)
    else:
       print 'Bad Domain'
\end{sageblock}

\noindent
Here is the code corresponding to the proof of Theorem~\ref{sp_area}:

\begin{sageblock}
#auto
@interact
def plotter(x = slider(pi/2,2*pi/3), y = slider(pi/3,2*pi/3)):
    if y <= x:        
        z = var('z')
        a, b, c = (x+y)/4, (x-y)/4, z/4
        expression(z) = tan(a+c)*tan(a-c)*tan(c+b)*tan(c-b)
        zmin, zmax = x, 2*pi/3
        color = (1,.25,0)
        plot(expression, z, zmin, zmax, rgbcolor=color)
            .show(xmin=pi/3,xmax=2*pi/3,ymin=0,ymax=.5)
    else:
       print 'Bad Domain'
\end{sageblock}
\bigskip

{\bf Acknowledgment.} We would like to thank Achill Sch\"{u}rmann and the anonymous referees for their helpful comments on the subject of this paper, and Andrew Ohana for helping us with some computations in Sage.
\bigskip

\bibliographystyle{plain}  
\bibliography{fukshansky_robins}    

\end{document}